\documentclass[a4,12pt]{article}
\usepackage{amsmath}
\usepackage{graphics}
\usepackage{latexsym}
\usepackage{amsfonts}
\numberwithin{equation}{section} \rightmargin 1.5cm \leftmargin
-1.5cm \textheight 23cm \textwidth 16cm \topmargin -1cm
\title{
Deterministic Minimax Impulse Control in Finite Horizon: the
Viscosity Solution Approach.}
\author{Brahim EL ASRI\\  Laboratoire Manceau de Math\'ematiques\\
Universit\'e du Maine\\ Avenue Olivier Messiaen\\ F-72085 Le Mans
Cedex 9 \\France\\ e-mail: brahim.el\_Asri@univ-lemans.fr }
\begin{document}
\date{}
\maketitle
\newtheorem{theo}{Theorem}
\newtheorem{problem}{Problem}
\newtheorem{pro}{Proposition}
\newtheorem{cor}{Corollary}
\newtheorem{axiom}{Definition}
\newtheorem{rem}{Remark}
\newtheorem{lem}{Lemma}
\newcommand{\brm}{\begin{rem}}
\newcommand{\erm}{\end{rem}}
\newcommand{\beth}{\begin{theo}}
\newcommand{\eeth}{\end{theo}}
\newcommand{\bl}{\begin{lem}}
\newcommand{\el}{\end{lem}}
\newcommand{\bp}{\begin{pro}}
\newcommand{\ep}{\end{pro}}
\newcommand{\bcor}{\begin{cor}}
\newcommand{\ecor}{\end{cor}}
\newcommand{\be}{\begin{equation}}
\newcommand{\ee}{\end{equation}}
\newcommand{\beq}{\begin{eqnarray*}}
\newcommand{\eeq}{\end{eqnarray*}}
\newcommand{\beqa}{\begin{eqnarray}}
\newcommand{\eeqa}{\end{eqnarray}}
\newcommand{\dg}{\displaystyle \delta}
\newcommand{\cm}{\cal M}
\newcommand{\cF}{{\cal F}}
\newcommand{\cR}{{\cal R}}
\newcommand{\bF}{{\bf F}}
\newcommand{\tg}{\displaystyle \theta}
\newcommand{\w}{\displaystyle \omega}
\newcommand{\W}{\displaystyle \Omega}
\newcommand{\vp}{\displaystyle \varphi}
\newcommand{\ig}[2]{\displaystyle \int_{#1}^{#2}}
\newcommand{\integ}[2]{\displaystyle \int_{#1}^{#2}}
\newcommand{\produit}[2]{\displaystyle \prod_{#1}^{#2}}
\newcommand{\somme}[2]{\displaystyle \sum_{#1}^{#2}}
\newlength{\inter}
\setlength{\inter}{\baselineskip} \setlength{\baselineskip}{7mm}
\newcommand{\no}{\noindent}
\newcommand{\rw}{\rightarrow}
\def \ind{1\!\!1}
\def \R{I\!\!R}
\def \N{I\!\!N}
\def \cadlag {{c\`adl\`ag}~}
\def \esssup {\mbox{ess sup}}
\begin{abstract}
We study here the impulse control minimax problem. We allow the cost
functionals and dynamics to be unbounded and hence the value
functions can possibly be unbounded. We prove that the value
function of the problem is continuous. Moreover, the value function
is characterized as the unique  viscosity solution of an Isaacs
quasi-variational inequality. This problem is in relation with an
application in mathematical finance.
\end{abstract}

\medskip

\no {$\bf Keywords$}: Impulse control; Robust control; Differential
games; Quasi-variational inequality; Viscosity solution

\section {Introduction}
In this paper we study an optimal impulse control problem with
finite horizon.

Optimal impulse control problems appear in many practical
situations. We refer the reader to \cite{[BL]} (and the references
cited therein) for extensive discussions. For deterministic
autonomous systems with infinite horizon, optimal impulse control
problems were studied in \cite{[B]}, and optimal control problems
with continuous, switching, and impulse controls were studied by the
author \cite{[Y11]} (see also \cite{[Y1]}). Differential games with
switching strategies in finite and infinite duration were also
studied \cite{[Y2], [Y3]}.  J. Yong, in [9], also studies
differential games where one person uses an impulse control and
other uses continuous controls. Recently El Farouq et al
\cite{[EBB]} extended the work of Yong \cite{[Y]} in the finite
horizon case but allowing general jumps. In all these works the
existence of the value functions of optimal impulse control problem
and uniqueness of viscosity solution are obtained assuming that the
dynamics and costs functionals are bounded and the impulse cost
function should not depend on y. \\

Our aim in this work is to relax the boundedness assumption on cost
functionals and the dynamics for impulse control problem and the
impulse cost function should depend on $y$.  Therefore the main
objective of our work, and this is the novelty of the paper, is to
show the existence of the value function and to characterize the
value function as the only solution in viscosity sense  of the
associated quasivariational inequality for the finite horizon
problem in suitable subclass of bounded from below continuous
functions, with linear growth when the dynamics unbounded and costs
functionals are bounded from below with linear growth and the
impulse cost function depends on $y$.

This paper is organized as follows:

In Section 2, we formulate the problem and we give the related
definitions.  In Section 3, we give some properties of the value
function, especially the dynamic programming principle. Further we
provide some estimates for the optimal strategy of the optimal
impulse control  problem which, in combination with the dynamic
programming principle, play a crucial role in the proof of the
existence of the value function. Section 4 is devoted to the
connection between the optimal impulse control  problem and
quasi-variational inequality. In Section 5, we show that the
solution of  QVIs is unique in the subclass of bounded from below
continuous functions which satisfy a linear growth condition. $\Box$
\section{Formulation of the problem and preliminary results}

\subsection{Setting of the problem}
Let a two-players differential game system be defined by the
solution of following dynamical equations
\begin{equation}\label{dynamic}\left\{
\begin{array}{l}
\dot{y}(t)=f(t,y(t),\tau(t))\\
y(t_0)=x\in \R^m,\\
y(t_k^+)=y(t_k^-)+g(t_k,y(t_k),\xi_k),\,\,t_k\geq t_0,\,\xi_k\neq 0,
\end{array}
\right.
\end{equation}
where $y(t)$ is the state of the system, with values in $\R^m$, at
time $t$, $x$ is the initial state. The time variable belongs to
$[t_0,T]$ where $0\leq t_0<T$.\\
The system is driven by two controls, a continuous control
$\tau(t)\in \mathcal{K}\in \R^m$, where $\mathcal{K}$ is compact
set, and an impulsive control defined by a double sequence $t_1, . .
. , t_k, . . . , \xi_1, . . . , \xi_k, . . .,k \in \N^*$, where
$t_k$ are the strategy, $t_k\leq t_{k+1}$ and $\xi_k\in\R^m$ the
control at time $t_k$ of the jumps in $y(t_k)$. Let
$(\delta,\xi):=((t_k)_{k\geq 1},(\xi_k)_{k\geq 1})$ the set of these
strategies denoted by
$\cal D$.\\
For any initial condition $(t_0, x)$ and controls $\tau(·)$ and
$\cal D$ generating a trajectory $y(·)$ of this system. The pay-off
is given by the following:

\begin{equation}\label{pay-off} J(t_0,x,\xi(.),\tau(.))=\integ{t_0}{T}\psi(s,y(s),\tau(s))ds
+\sum_{k\geq 1} C(t_k,y(t_k),\xi_k)\ind_{[t_{k}\leq T]}+G(y(T)).
\end{equation}
The term $C(t_k,y(t_k),\xi_k)$ is called the impulse cost. It is the
cost when player-$\xi$ makes an impulse $\xi_k$ at time $t_k$. In
the game, player-$\xi$ would like to minimize the pay-off
(\ref{pay-off}) by choosing suitable impulse control $\xi(.)$,
whereas player-$\tau$ wants to maximize the pay-off (\ref{pay-off})
by choosing a proper control $$\tau(.)\in \Omega=\{\mbox{measurable
functions}\quad [t_0,T]\rightarrow \mathcal{K}\}.$$ We shall
sometimes write $\tau\in \Omega$ instead of $\tau(.)\in \Omega$.\\

We now define the admissible strategies $\varphi$ for the minimizing
impulse control $\cal D$, as non-anticipative strategies. We shall
let ${\cal D}_a$ be the set of all such non-anticipative strategies.
\begin{axiom}
A map $\varphi: \Omega\rightarrow {\cal D}_a$ is called a
non-anticipative strategy if for any two controls $\tau_1(.)$ and
$\tau_2(.)$, and any $t\in[t_0,T ]$, the condition on their
restrictions to $[t_0, t[: \tau_1|_{[t_0,t [} = \tau_2|_{[t_0,t [} $
implies $\varphi(\tau_1)|_{[t_0,t ]} = \varphi(\tau_2)|_{[t_0,t ]}$.
\end{axiom}

Accordingly, we define the value function of the problem
$v:[0,T]\times \R^m\rightarrow \R$ as
$$v(t_0,x)=\inf\limits_{\varphi\in {\cal D}_a}\sup\limits_{\tau(.)\in \Omega}J(t_0,x,\xi(.),\tau(.))$$

\subsection{Assumptions} Throughout
this paper $T$ (resp. $m$) is a fixed real (resp. integer) positive
constant. Let us now consider the followings:
\medskip

\indent $(1)$ $f:[0,T]\times\R^m\times \mathcal{K}\rightarrow
\R^{m}$ and $g:[0,T]\times\R^m \times \R^m\rightarrow \R^{m}$ are
two continuous functions for which there exists a constant $C\geq 0$
such that for any $t\in [0, T]$, $\tau \in \mathcal{K}$ and $\xi, x,
x'\in \R^m$\be \label{regbs1}|f(t,x,\tau)|+ |g(t,x,\xi)|\leq
C(1+|x|) \quad \mbox{ and } \quad
|g(t,x,\xi)-g(t,x',\xi)|+|f(t,x,\tau)-f(t,x',\tau)|\leq C|x-x'|\ee

$(2)$ $C:[0,T]\times\R^m \times \R^m\rightarrow \R^{}$, is
continuous with respect to $t$ and $\xi$ uniformly in $y$ with
linear growth $$|C(t,x,\xi)|\leq C(1+|x|),\, \, \forall (t,x,\xi)\in
[0,T]\times \R^m\times  \R^m.$$ Moreover we assume that there exists
a constant $\alpha
>0$ such that for any $(t,x,\xi)\in [0,T]\times \R^m\times \R^m$,
\be C(t,x,\xi)\geq \alpha.\ee

$(3)$ $\psi:[0,T]\times\R^m\times \mathcal{K}\rightarrow \R$ is
continuous  with respect to $t$ and $\tau$ uniformly in $y$ with
linear growth, \be \label{polycond} |\psi(t,x,\tau)|\leq C(1+|x|),\,
\, \forall (t,x,\tau)\in [0,T]\times \R^m\times  \mathcal{K},\ee and
is bounded from below.

$(4)$$G:\R^m\rightarrow \R$ is uniformly continuous  with linear
growth \be \label{polycond1}|G(x)|\leq C(1+|x|),\, \, \forall x\in
 \R^m, \ee
 and
is bounded from below.
%

These properties of $f$ and $g$ imply in particular that $y(t)_{0\le
t\leq T}$ solution of the standard DE (\ref{dynamic}) exists and is
unique, for any $t\in [0, T]$ and $x\in \R^m$.
\medskip

\subsection{Admissible strategies}
We want to investigate the problem of minimizing
$\sup\limits_{\tau\in \Omega}J$ through the impulse control. We mean
to allow closed loop strategies for the minimizing control. We
remark that, being only interested in the $\inf \sup$ problem, and
not a possible saddle point.
\begin{theo}
Under the standing assumptions (Sect. 2.2) the value function $v$ is
bounded from below with linear growth.
\end{theo}
$Proof$: Consider the particular strategy in ${\cal D}_a$ is the one
where we have no impulse time. In this case $\sum_{k\geq 1}
C(t_k,y(t_k),\xi_k)\ind_{[t_{k}\leq T]}=0.$

$$ v(t,x)\leq \sup\limits_{\tau \in \Omega}[\integ{t}{T}\psi(s,y(s),\tau(s))ds +G(y(T))].
$$

Since $\psi$ and $G$ are linear growth, then
$$ v(t,x)\leq \integ{t}{T}C(1+|y(s)|)ds +C(1+|y(T)|).
$$
Now by using standard estimates from ODE, Gronwall's Lemma and the
strategy where we have no impulse time, we can show that
$$|y(t)|\leq C(1+|x|),$$
where $C$ is constant of $T$. Hence using this estimate we get

$$ v(t,x)\leq C(1+|x|).$$
 On the other hand, since the cost $C(t_k,y(t_k),\xi_k)$
are non negative functions and since $\psi$ and $G$ are bounded from
below,
then  $v$ is bounded from below.$\Box$\\

 We are now giving some properties of the admissible strategy.

 \begin{pro}
\label{optimal-s} Let $(\delta,\xi)=((t_n)_{n\geq 1},(\xi_n)_{n\geq
1})$ be an optimal strategy. Then:
$$\sum_{k\geq 1} C(t_k,y(t_k),\xi_k)\ind_{[t_{k}\leq T]}\leq C(1+|x|).$$

There exists a positive constant $C$  which dose not depend on $t$
and $x$ such that: \be \label{estiopt1}\forall\, n\geq
1,\,\,\ind_{[t_{n}\leq T]}\leq \frac {C(1+|x|)}{n}.\ee \ep $Proof$:
Recall the characterization of (\ref{pay-off}) that reads as:
 $$v(t,x)=\inf\limits_{\varphi\in {\cal D}_a}\sup\limits_{\tau(.)\in \Omega}[\integ{t}{T}\psi(s,y(s),\tau(s))ds
+\sum_{k\geq 1} C(t_k,y(t_k),\xi_k)\ind_{[t_{k}\leq T]}+G(y(T))].$$
Now if $((t_n)_{n\geq 1},(\xi_n)_{n\geq 1})$ is the optimal strategy
then we have:
$$\integ{t}{T}\psi(s,y(s),\tau(s))ds
+\sum_{k\geq 1} C(t_k,y(t_k),\xi_k)\ind_{[t_{k}\leq T]}+G(y(T))\leq
v(t,x).$$

Since $v(t,x)\leq C(1+|x|)$ and since $\psi$ and $G$ are bounded
from below then we have

$$\sum_{k\geq 1} C(t_k,y(t_k),\xi_k)\ind_{[t_{k}\leq T]}\leq
C(1+|x|).$$

Next we show (\ref{estiopt1}). Taking into account that
$C(t,y,\xi)\geq \alpha
>0$ for any $(t,y,\xi)\in[t_0,T]\times\R^m\times\R^m$
 we obtain:
$$\sum_{k\geq 1} \alpha \ind_{[t_{k} \leq T]}\leq
C(1+|x|).$$

But for any $k\le n$, $[t_{n}\leq T]\subset[t_{k}\leq T]$
then:$$\begin{array}{l} \alpha n \ind_{[t_{n}\leq T]}\leq C(1+|x|).
\end{array}
$$
Finally taking into account $\alpha>0$, we obtain the desired
result.$\Box$

It may be to the best advantage of the minimizer to make a jump at
some time $t$, immediately followed, at the same time, by another
jump, and so on. As any such jump entails a cost not less than
$\alpha$, from Proposition 1 
the number of jumps may be restricted, with no loss of generality,
to be less than $\frac{C}{\alpha}(1+|x|)$. To allow for the
possibility of several successive but simultaneous jumps, we proceed
as follows. Let
$$E=\left\{\xi\in \bigcup_{n=1}(\R^m)^n,\quad\mbox{such that}\quad \forall n, \, \ind_{[t_{n}\leq T]}\leq \frac {C(1+|x|)}{n}\right\}.$$

We also state the following definition:
\begin{axiom}
For any function $v:[t_0,T]\times \R^m\rightarrow \R$, let the
operator $N$ be given by
$$N[v](t,x)=\inf\limits_{\xi\in E}[v(t,x+g(t,x,\xi))+C(t,x,\xi)].$$
\end{axiom}

\section{The value function}
\subsection{Dynamic programming principle}
The dynamic programming principle is a well-known property in
 optimal  impulse control. In our optimal control problem, it is
formulated as follows:
\begin{theo}(\cite{[EBB]}, Proposition 3.1)\label{principle-dyn}
The value function $v(.,.)$ satisfies the following optimality
principle:\\
for all $t\leq t'\, \in [t_0,T[$ and $x\in \R^m$,
$$v(t,x)=\inf\limits_{\varphi\in {\cal D}_a}\sup\limits_{\tau\in \Omega}[\integ{t}{t'}\psi(s,y(s),\tau(s))ds
+\sum_{k\geq 1,\,t_k<t'} C(t_k,y(t_k),\xi_k)\ind_{[t_{k}\leq
T]}+\ind_{[t' \leq T]}v(t',y(t'))],$$

and
$$v(t,x)=\inf\limits_{\varphi\in {\cal D}_a}\sup\limits_{\tau\in \Omega}[\integ{t}{t_n}\psi(s,y(s),\tau(s))ds
+\sum_{1\leq k<n} C(t_k,y(t_k),\xi_k)\ind_{[t_{k}\leq T]}+\ind_{[t_n
\leq T]}v(t_n,y(t_n))],$$ where $((t_n)_{n\geq 1},(\xi_n)_{n\geq
1})$ be an admissible control.
\end{theo}

\begin{pro}
The value function v(.,.) has the following property: \\
for all $t\in [t_0,T]$ and $x\in \R^m$,
$$v(t,x)\leq N[v](t,x).$$
\end{pro}
 $Proof:$
Assume first that for some $x$ and $t$:
$$v(t,x)> N[v](t,x).$$
Then we have for $t\leq t'$:\begin{eqnarray*} &&
\inf\limits_{\varphi\in {\cal D}_a}\sup\limits_{\tau\in
\Omega}[\integ{t}{t'}\psi(s,y(s),\tau(s))ds +\sum_{k\geq 1,\,t_k<t'}
C(t_k,y(t_k),\xi_k)\ind_{[t_{k}\leq T]}+\ind_{[t'\leq T]}v(t',y(t'))]\\
&& \quad >\inf\limits_{\xi\in E}[v(t,x+g(t,x,\xi))+C(t,x,\xi)].
\end{eqnarray*}
Among the admissible strategy $\varphi^\epsilon$'s  there are those
that place a jump at time t.
\begin{eqnarray*}
&& \sup\limits_{\tau\in \Omega}[\integ{t}{t'}\psi(s,y(s),\tau(s))ds
+\sum_{k\geq 1,\,t_k<t'}
C(t_k,y(t_k),\xi_k)\ind_{[t_{k}<T]}+\ind_{[t'<T]}v(t',y(t'))]\\
&& \quad > v(t,x+g(t,x,\xi))+C(t,x,\xi)-\epsilon.
\end{eqnarray*}
Now, pick $\tau_1$ such that
\begin{eqnarray*}
&& \integ{t}{t'}\psi(s,y(s),\tau_1(s))ds +\sum_{k\geq 1,\,t_k<t'}
C(t_k,y(t_k),\xi_k)\ind_{[t_{k}\leq T]}+\ind_{[t'\leq T]}v(t',y(t')) +\epsilon\\
&& \quad \geq \sup\limits_{\tau\in
\Omega}[\integ{t}{t'}\psi(s,y(s),\tau(s))ds +\sum_{k\geq 1,\,t_k<t'}
C(t_k,y(t_k),\xi_k)\ind_{[t_{k}\leq T]}+\ind_{[t'\leq
T]}v(t',y(t'))],
\end{eqnarray*}
which implies that:
\begin{eqnarray*}
&& \integ{t}{t'}\psi(s,y(s),\tau_1(s))ds +\sum_{k\geq 1,\,t_k<t'}
C(t_k,y(t_k),\xi_k)\ind_{[t_{k}\leq T]}+\ind_{[t'\leq T]}v(t',y(t')) +\epsilon\\
&& \quad > v(t,x+g(t,x,\xi))+C(t,x,\xi)-\epsilon.
\end{eqnarray*}
Choosing now $t'=t$, yields the relation
$$ \epsilon+v(t,x+g(t,x,\xi))> v(t,x+g(t,x,\xi))+C(t,x,\xi)-\epsilon.$$
By sending $\epsilon\rightarrow0$, we obtain $C(t,x,\xi)<0 $, which
is a contradiction.$\Box$
\subsection{Continuity of value function}
In this section we prove the continuity of the value function. The
main result of this section can be stated as follows.\\

We first present some preliminary results on $y(.).$ Consider now
the control ${\cal D}_a$, composed of jumps instants
$t_1,t_2,...,t_n$ in the interval $[t,T]$, with jumps
$\xi_1,\xi_2,...,\xi_n,$ and let $y_1(.)$ and $y_2(.)$ be the
trajectories generated by ${\cal D}_a$, from $y_i(t)=x_i,\,i=1,2.$
\begin{lem}
There exists a constant $C$ such that for any $s\in [t,T]$,
$x_1,x_2\in \R^m$, and $k\in\{1,2...,n\}$
\begin{equation}\label{estimat2}
|y_1(s)-y_2(s)|\leq \exp(C(s-t))(1+C)^n|x_1-x_2|. \Box
\end{equation}
\end{lem}
$Proof:$ By the Lipschitz continuity of $f$ and Gronwall's Lemma, we
have

$$|y_1(s)-y_2(s)|\leq \exp(C(s-t))(1+C)|x_1-x_2|,\quad \forall s\in[t,t_1].$$
Next let us show for an impulse time
$$|y_1(t_k^+)-y_2(t_k^+)|\leq \exp(C(t_k-t))(1+C)^k|x_1-x_2|.$$
Looking more carefully at the first jump and using the Lipschitz
continuity of $g$, we have
\begin{equation}
\begin{array}{lll}
\label{viscder}
|y_1(t_1^+)-y_2(t_1^+)|&=|y_1(t_1^-)+g(t_1^-,y_1(t_1^-),\xi_1)-y_2(t_1^-)+g(t_1^-,y_2(t_1^-),\xi_1)|\\&\leq
(1+C)|y_1(t_1^-)-y_2(t_1^-)|
\\&\leq \exp(C(t_1-t))(1+C)|x_1-x_2|.
\end{array}
\end{equation}
The above assertion is obviously true for $k = 1$. Suppose now it
holds true at step $k$. Then, at step $k + 1$,
\begin{equation}
\begin{array}{lll}
\label{viscder} |y_1(t_{k+1}^+)-y_2(t_{k+1}^+)|&\leq
(1+C)|y_1(t_{k+1}^-)-y_2(t_{k+1}^-)|\\&\leq
(1+C)|y_1(t_k^+)-y_2(t_k^+)|\exp(C(t_{k+1}^- -t_k^+))
\\&\leq  \exp(C(t_{k+1}-t))(1+C)^{k+1}|x_1-x_2|.
\end{array}
\end{equation}
Finally $$|y_1(s)-y_2(s)|\leq
\exp(C(s-t))(1+C)^n|x_1-x_2|,\quad\forall s\in [0,T]. \Box$$

\medskip

We are now ready to give the main Theorem of this article. \beth The
value function $v:[0,T]\times \R^m\rightarrow \R$ is continuous in
$t$ and $x$.\eeth $Proof$: Let us consider $\epsilon >0$ and
$(t',x')\in B((t,x),\epsilon)$ and let us consider the following set
of strategies:$$ \tilde D_a:=\left\{(\delta,\xi)=((t_n)_{n\geq 1},
(\xi_n)_{n\geq 0}) \in {\cal D}_a \mbox{ such that } \forall n\geq
1, \ind_{[\tau_n\leq T]} \leq \frac{C(1+(\epsilon
+|x|))}{n}\right\}.$$ The strategy optimal $(\delta,\xi)$
belongs to $\tilde D_a$ from Proposition 1.\\
 First let us show that $v$ is upper semi-continuous. Recall
the characterization of dynamical programming principle that reads
as
$$v(t,x)=\inf\limits_{\varphi\in {\tilde  D}_a}\sup\limits_{\tau\in \Omega}[\integ{t}{t_n}\psi(s,y(s),\tau(s))ds
+\sum_{1\leq k<n} C(t_k,y(t_k),\xi_k)\ind_{[t_{k}\leq T]}+\ind_{[t_n
\leq T]}v(t_n,y(t_n))],$$
$$v(t',x')=\inf\limits_{\varphi\in {\tilde  D}_a}\sup\limits_{\tau\in \Omega}[\integ{t'}{t_n}\psi(s,y'(s),\tau(s))ds
+\sum_{1\leq k<n} C(t_k,y'(t_k),\xi_k)\ind_{[t_{k}\leq
T]}+\ind_{[t_n \leq T]}v(t_n,y'(t_n))].$$ Fix an arbitrary
$\epsilon^1>0$. Let $\varphi=((t_n)_{n\geq 1},(\xi_n)_{n\geq 1})$
belongs to $\tilde D_a$ such that
\begin{eqnarray*}
&&\sup\limits_{\tau\in \Omega}[\integ{t}{t_n}\psi(s,y(s),\tau(s))ds
+\sum_{1\leq k<n} C(t_k,y(t_k),\xi_k)\ind_{[t_{k}\leq
T]}+\ind_{[t_n \leq T]}v(t_n,y(t_n))] \\
&& \quad \leq \inf\limits_{\varphi\in {\tilde
D}_a}\sup\limits_{\tau\in
\Omega}[\integ{t}{t_n}\psi(s,y(s),\tau(s))ds +\sum_{1\leq k<n}
C(t_k,y(t_k),\xi_k)\ind_{[t_{k}\leq T]}+\ind_{[t_n \leq
T]}v(t_n,y(t_n))]+\epsilon^1\\ && \quad =v(t,x)+\epsilon^1.
\end{eqnarray*}
Also,
$$v(t',x')\leq \sup\limits_{\tau\in
\Omega}[\integ{t'}{t_n}\psi(s,y'(s),\tau(s))\ind_{[s\geq t']}ds
+\sum_{1\leq k<n} C(t_k,y'(t_k),\xi_k)\ind_{[t_{k}\leq
T]}+\ind_{[t_n \leq T]}v(t_n,y'(t_n))].$$ Now pick $\tau_1$ such
that

\begin{eqnarray*}
&& \sup\limits_{\tau\in
\Omega}[\integ{t'}{t_n}\psi(s,y'(s),\tau(s))\ind_{[s\geq t']}ds
+\sum_{1\leq k<n} C(t_k,y'(t_k),\xi_k)\ind_{[t_{k}\leq
T]}+\ind_{[t_n \leq
T]}v(t_n,y'(t_n))]\\
&& \quad \leq \integ{t'}{t_n}\psi(s,y'(s),\tau^1(s))\ind_{[s\geq
t']}ds +\sum_{1\leq k<n} C(t_k,y'(t_k),\xi_k)\ind_{[t_{k}\leq
T]}+\ind_{[t_n \leq T]}v(t_n,y'(t_n))+\epsilon^1.
\end{eqnarray*}

Then
\begin{eqnarray*}
&& v(t',x')-v(t,x)\leq
\integ{t'}{t_n}\psi(s,y'(s),\tau^1(s))\ind_{[s\geq t']}ds
+\sum_{1\leq k<n} C(t_k,y'(t_k),\xi_k)\ind_{[t_{k}\leq T]}\\&&
\qquad\qquad \qquad+\ind_{[t_n \leq T]}v(t_n,y'(t_n))]-
\integ{t}{t_n}\psi(s,y(s),\tau^1(s))ds -\sum_{1\leq k<n}
C(t_k,y(t_k),\xi_k)\ind_{[t_{k}\leq T]}\\&& \qquad\qquad
\qquad-\ind_{[t_n \leq T]}v(t_n,y(t_n))+2\epsilon^1.
\end{eqnarray*}
Next w.l.o.g we assume that $t'<t$. Then we deduce that:
\begin{equation}\label{con_sup}
\begin{array}{ll}v(t',x')-v(t,x)&\leq
\integ{t_0}{t_n}\displaystyle\{(\psi(s,y'(s),\tau^1(s))-\psi(s,y(s),\tau^1(s)))\ind_{[s\geq t]}+\psi(s,y'(s),\tau^1(s))\ind_{[t'\leq s< t]}\}ds\\
{}&\qquad +\displaystyle\sum_{1\leq k<n}\{
C(t_k,y'(t_k),\xi_k)-C(t_k,y(t_k),\xi_k)\}\ind_{[t_{k}\leq
T]}\\{}&\qquad +\displaystyle\ind_{[t_n \leq
T]}\{v(t_n,y'(t_n))-v(t_n,y(t_n))\}+2\epsilon^1\\
{}&\leq \integ{t_0}{t_n}\displaystyle\{|\psi(s,y'(s),\tau^1(s))-\psi(s,y(s),\tau^1(s))|\ind_{[s\geq t]}+|\psi(s,y'(s),\tau^1(s))|\ind_{[t'\leq s< t]}\}ds\\
{}&\qquad +\displaystyle n\max\limits_{1\leq k\leq n}
|C(t_k,y'(t_k),\xi_k)-C(t_k,y(t_k),\xi_k)|\\{}&\qquad
+\displaystyle\ind_{[t_n \leq
T]}\{|v(t_n,y'(t_n))|+|v(t_n,y(t_n))|\}+2\epsilon^1.
\end{array}
\end{equation}
Using the uniform continuity of $\psi$, $C$ in $y$ and property
(\ref{estimat2}), then the right-hand side of (\ref{con_sup}), the
first and the second term converges to 0 as $t'\rightarrow t$ and
$x'\rightarrow x$.\\
Now let us focus on the last one. Since $(\delta,\xi) \in \tilde
D_a$  then
$$\ind_{[t_n
\leq T]}\{|v(t_n,y'(t_n))|+|v(t_n,y(t_n))|\}\leq
C\frac{(1+|x|^2+|x'|^2)}{n},$$

\noindent where $C$ is a constant which comes from the linear growth
of $\psi$ and $G$. Taking the limit as $(t',x')\rw (t,x)$ we obtain:
$$
\limsup_{(t',x')\rw (t,x)}v(t',x')\leq v(t,x)
+C\frac{(1+|x|^2+|x'|^2)}{n}+2\epsilon^1.$$As $n$ and $\epsilon^1$
are arbitrary then sending $n\rw +\infty$ and $\epsilon^1\rightarrow
0,$ to obtain:
$$
\limsup_{(t',x')\rw (t,x)}v(t',x')\leq v(t,x).$$ Therefore $v$ is
upper semi-continuous.\\
Now we show that $v$ is lower semi-continuous.\\
 Fix an arbitrary
$\epsilon_2>0$. Let $\varphi_2=((t_n)_{n\geq 1},(\xi_n)_{n\geq 1})$
belongs to $\tilde D_a$ such that
\begin{eqnarray*}
&&\sup\limits_{\tau\in
\Omega}[\integ{t'}{t_n}\psi(s,y'(s),\tau(s))ds +\sum_{1\leq k<n}
C(t_k,y'(t_k),\xi_k)\ind_{[t_{k}\leq
T]}+\ind_{[t_n \leq T]}v(t_n,y'(t_n))] \\
&& \quad \leq \inf\limits_{\varphi_2\in {\tilde
D}_a}\sup\limits_{\tau\in
\Omega}[\integ{t}{t_n}\psi(s,y'(s),\tau(s))ds +\sum_{1\leq k<n}
C(t_k,y'(t_k),\xi_k)\ind_{[t_{k}\leq T]}+\ind_{[t_n \leq
T]}v(t_n,y(t_n))]+\epsilon^1\\ && \quad =v(t',x')+\epsilon_2.
\end{eqnarray*}
Also, $$v(t,x)\leq \sup\limits_{\tau\in
\Omega}[\integ{t}{t_n}\psi(s,y(s),\tau(s))\ind_{[s\geq t]}ds
+\sum_{1\leq k<n} C(t_k,y(t_k),\xi_k)\ind_{[t_{k}\leq T]}+\ind_{[t_n
\leq T]}v(t_n,y(t_n))],$$ now, pick $\tau_2$ such that

\begin{eqnarray*}
&& \sup\limits_{\tau\in
\Omega}[\integ{t}{t_n}\psi(s,y(s),\tau(s))\ind_{[s\geq t]}ds
+\sum_{1\leq k<n} C(t_k,y(t_k),\xi_k)\ind_{[t_{k}\leq T]}+\ind_{[t_n
\leq
T]}v(t_n,y(t_n))]\\
&& \quad \leq \integ{t}{t_n}\psi(s,y(s),\tau_2(s))\ind_{[s\geq t]}ds
+\sum_{1\leq k<n} C(t_k,y(t_k),\xi_k)\ind_{[t_{k}\leq T]}+\ind_{[t_n
\leq T]}v(t_n,y(t_n))-\epsilon_2.
\end{eqnarray*}

Then
\begin{eqnarray*}
&& v(t',x')-v(t,x)\geq \integ{t'}{t_n}\psi(s,y'(s),\tau_2(s))ds
+\sum_{1\leq k<n} C(t_k,y'(t_k),\xi_k)\ind_{[t_{k}\leq
T]}+\ind_{[t_n \leq
T]}v(t_n,y'(t_n))]\\
&& \qquad\qquad \qquad- \integ{t}{t_n}\psi(s,y(s),\tau_2(s))ds
-\sum_{1\leq k<n} C(t_k,y(t_k),\xi_k)\ind_{[t_{k}\leq T]}-\ind_{[t_n
\leq T]}v(t_n,y(t_n))-2\epsilon_2.
\end{eqnarray*}
Next w.l.o.g we assume that $t'<t$. Then we deduce that:
\begin{equation}\label{con_sup1}
\begin{array}{ll}v(t',x')-v(t,x)&\geq
\integ{t_0}{t_n}\displaystyle\{(\psi(s,y'(s),\tau_2(s))-\psi(s,y(s),\tau_2(s)))\ind_{[s\geq t]}+\psi(s,y'(s),\tau_2(s))\ind_{[t'\leq s< t]}\}ds\\
{}&\qquad +\displaystyle\sum_{1\leq k<n}\{
C(t_k,y'(t_k),\xi_k)-C(t_k,y(t_k),\xi_k)\}\ind_{[t_{k}\leq
T]}\\{}&\qquad +\displaystyle\ind_{[t_n \leq
T]}\{v(t_n,y'(t_n))-v(t_n,y(t_n))\}-2\epsilon_2\\
{}&\geq -\integ{t_0}{t_n}\displaystyle\{|\psi(s,y'(s),\tau_2(s))-\psi(s,y(s),\tau_2(s))|\ind_{[s\geq t]}+|\psi(s,y'(s),\tau_2(s))|\ind_{[t'\leq s< t]}\}ds\\
{}&\qquad -\displaystyle n\max\limits_{1\leq k\leq n}
|C(t_k,y'(t_k),\xi_k)-C(t_k,y(t_k),\xi_k)|\\{}&\qquad -\displaystyle
\ind_{[t_n \leq T]}\{|v(t_n,y'(t_n))|+|v(t_n,y(t_n))|\}-2\epsilon_2.
\end{array}
\end{equation}
Using the uniform continuity of $\psi$, $C$ in $y$ and property
(\ref{estimat2}). Then the right-hand side of (\ref{con_sup}) the
first and the second term converges to 0 as $t'\rightarrow t$ and
$x'\rightarrow x$.\\
Now let us focus on the last one. Since $((t_n)_{n\geq
1},(\xi_n)_{n\geq 1})$ be a admissible control then
$$-\ind_{[t_n
\leq T]}\{|v(t_n,y'(t_n))|+|v(t_n,y(t_n))|\}\geq -
C\frac{(1+|x|^2+|x'|^2)}{n},$$

where $C$ is a constant which come from the linear growth of $\psi$
and $G$, taking the limit as $(t',x')\rw (t,x)$ to obtain:
$$
\liminf_{(t',x')\rw (t,x)}v(t',x')\geq v(t,x)
-C\frac{(1+|x|^2+|x'|^2)}{n}-2\epsilon_2.$$As $n$ and $\epsilon_2$
are arbitraries then putting $n\rw +\infty$ and
$\epsilon_2\rightarrow 0$ to obtain:
$$
\liminf_{(t',x')\rw (t,x)}v(t',x')\geq v(t,x).$$ Therefore $v$ is
lower semi-continuous. We then proved that $v$ is  continuous.$\Box$\\

\subsection{Terminal value}
Because of the possible jumps at the terminal time T, it is easy to
see that, in general, $v(t,x)$ does not tend to $G(x)$ as t tends to
T. Extend the set of jumps to include jumps of zero, meaning no
jump. Call this extended set $E_0$, extend trivially the operator
$N$ to a function independent from $t$, and let

\begin{equation}\label{Valeur-ter} G_1(x)=\inf\limits_{\xi \in
E_0}[G(x+g(T,x,\xi))+C(T,x,\xi)]=\min\{G(x),N[G](T,x)\}.
\end{equation}
We know that $G$ and $C$ are uniformly continuous in $x$ then
$G_1(x)$ is continuous. We claim
\begin{lem}

$$v(t,x)\rightarrow G_1(x)\quad \mbox{as} \quad t\rightarrow T.$$
\end{lem}
$Proof$: Fix $(t,x)$ and a strategy $\varphi$. As in the previous
proof, for each $\tau(·)$, gather all jumps of $\varphi(\tau)$ if
any, in jump $\xi_1$ at the time T. Then we have
$$|J(t,x,\varphi,\tau)-G(x+g(T,x,\xi_1))-C(T,x,\xi_1)|\leq C_x(T-t)$$
or
$$J(t,x,\varphi,\tau)=G(x+g(T,x,\xi_1))+C(T,x,\xi_1)+O(T-t).$$
The right hand side above only depends on $\xi_1$, not on $\tau(.)$
itself. It follows that
$$
\begin{array}{ll}\inf\limits_{\varphi}\sup\limits_{\tau}J(t,x,\varphi,\tau)&=\inf\limits_{\xi\in
E_0}[G(x+g(T,x,\xi))+C(T,x,\xi)]+O(T-t)\\ {}&=G_1(x)+O(T-t).
\end{array}$$
The result follows letting $t\rightarrow T.\Box$
\section{Viscosity characterization of the value function}

In this section we prove that the value function $v$ is a viscosity
solution of the Hamilton-Jacobi-Isaacs quasi-variational inequality,
that we replace by an equivalent QVI easier to investigate.\\

We now consider the following quasi-variational inequality (Isaacs
equation): \be
\begin{array}{l}\label{Issacs}
\displaystyle \max\left\{\min\limits_{\tau\in{\cal
K}}\left[-\frac{\partial v}{\partial t}- \frac{\partial v}{\partial
x}f(t,x,\tau)-\psi(t,x,\tau)\right],\right.\\\qquad\qquad \qquad
\left.v(t,x)-N[v](t,x)\right\}= 0,
\end{array}
\ee with the terminal condition: $v(t,x)=G_1(x),\, x\in \R^m$, where
$G_1$ is given by (\ref{Valeur-ter}).\\
Notice that it follows from hypothesis that the term in square
brackets in (\ref{Issacs}) above is continuous with respect to
$\tau$ so that the minimum in $\tau$ over the compact $\cal K$
exists.\\
Recall the notion of viscosity solution of QVI (\ref{Issacs}).
\begin{axiom} Let $v$ be a  continuous function defined on
$[0,T]\times \R^m$, $\R$-valued and such that $v(T,x)=G_1(x)$ for
any $x\in \R^m$. The $v$ is called:
\begin{itemize}
\item [$(i)$] A viscosity supersolution  of (\ref{Issacs})
if for any $(\overline{t},\overline{x})\in [t_0,T[\times \R^m$ and
any function $\varphi \in C^{1,2}([t_0,T[\times \R^m)$ such that
$\varphi(\overline{t},\overline{x})=v(\overline{t},\overline{x})$
and $(\overline{t},\overline{x})$ is a local maximum of $\varphi
-v$, we have: \be
\begin{array}{l}\label{Issacs1}
\displaystyle \max\left\{\min\limits_{\tau\in{\cal
K}}\left[-\frac{\partial v}{\partial t}- \frac{\partial v}{\partial
x}f(\overline{t},\overline{x},\tau)-\psi(\overline{t},\overline{x},\tau)\right],\right.\\\qquad\qquad
\qquad
\left.v(\overline{t},\overline{x})-N[v](\overline{t},\overline{x})\right\}\geq
0.
\end{array}
\ee

\item [$(ii)$] A viscosity  subsolution of (\ref{Issacs})
if for any $(\overline{t},\overline{x})\in [t_0,T[\times \R^m$ and
any function $\varphi \in C^{1,2}([t_0,T[\times \R^m)$ such that
$\varphi(\overline{t},\overline{x})=v(\overline{t},\overline{x})$
and $(\overline{t},\overline{x})$ is a local minimum of $\varphi
-v$, we have: \be
\begin{array}{l}\label{Issacs2}
\displaystyle \max\left\{\min\limits_{\tau\in{\cal
K}}\left[-\frac{\partial v}{\partial t}- \frac{\partial v}{\partial
x}f(\overline{t},\overline{x},\tau)-\psi(\overline{t},\overline{x},\tau)\right],\right.\\\qquad\qquad
\qquad
\left.v(\overline{t},\overline{x})-N[v](\overline{t},\overline{x})\right\}\leq0.
\end{array}
\ee
\item [$(iii)$] A viscosity solution if it is both a viscosity supersolution and
subsolution. $\Box$
\end{itemize}
\end{axiom}
\begin{theo}
The function: $(t,x)\rightarrow v(t,x)$ is viscosity solution of the
quasi-variational inequality (\ref{Issacs}).
\end{theo}
$Proof$: The viscosity property follows from the dynamic programming
principle and is proved in \cite{[EBB]}.$\Box$\\

Now we give an equivalent of quasi-variational inequality
(\ref{Issacs}). In this section, we consider the new function
$\Gamma$ given by the classical change of variable $\Gamma(t,x) =
\exp(t)v(t, x)$, for any $t\in[t_0,T ]$ and $x\in \R^m$. Of course,
the function $\Gamma$ is bounded from below and continuous with
respect to its arguments.\\ A second property is given by the

\begin{pro}
$v$ is a viscosity solution of (\ref{Issacs}) if and only if
$\Gamma$ is a viscosity solution to the following quasi-variational
inequality in $[t_0,T [\times \R^m$, \be
\begin{array}{l}\label{Issacs3}
\displaystyle \max\left\{\min\limits_{\tau}\left[-\frac{\partial
\Gamma}{\partial t}+\Gamma(t,x)- \frac{\partial \Gamma}{\partial
x}f(t,x,\tau)-\exp(t)\psi(t,x,\tau)\right],\right.\\\qquad\qquad
\qquad \left.\Gamma(t,x)-M[\Gamma](t,x)\right\}= 0,
\end{array}
\ee where $M[\Gamma](t,x)=\inf\limits_{\xi\in
E}[\Gamma(t,x+g(t,x,\xi))+\exp(t)C(t,x,\xi)].$ The terminal
condition for $\Gamma$ is: $\Gamma(T,x)=\exp(T)G_1(x)$ in
$\R^m.$$\Box$
\end{pro}

\section{Uniqueness of the solution of quasi-variational inequality} We are going now to address the
question of uniqueness of the viscosity solution of
quasi-variational inequality (\ref{Issacs}). We have the following:

\beth \label{uni1}The solution in viscosity sense of
quasi-variational inequality (\ref{Issacs}) is unique in the space
of continuous functions on $[t_0,T]\times R^m$ which satisfy a
linear growth condition, i.e., in the space
$$\begin{array}{l}{\cal C}:=\{\varphi: [0,T]\times \R^m\rightarrow
\R, \mbox{ continuous and for any }\\\qquad \qquad\qquad(t,x), \,
\varphi(t,x)\leq C(1+|x|) \mbox{ for some constants } C
\quad\mbox{and bounded from below} \}.\end{array}$$\eeth {\it
Proof}. We will show by contradiction that if $u$ and $w$ is a
subsolution and a supersolution respectively for (\ref{Issacs3})
then  $u\leq w$. Therefore if we have two solutions of
(\ref{Issacs3}) then they are obviously equal. Actually for some
$R>0$ suppose there exists
$(\overline{t},\overline{x})\in[t_0,T]\times B_R\times {\cal I}$
$(B_R := \{x\in \R^m; |x|<R\})$ such that:
\begin{equation}
\label{comp-equ1}
\max\limits_{t,x}(u(t,x)-w(t,x))=u(\overline{t},\overline{x})-w(\overline{t},\overline{x})=\eta>0.
\end{equation}Let us take
$\theta$, $\lambda$ and $\beta \in (0,1]$ small enough. 
 Then, for a small $\epsilon>0$, let us define:
\begin{equation}
\label{phi}
\Phi_{\epsilon}(t,x,y)=(1-\lambda)u(t,x)-w(t,y)-\frac{1}{2\epsilon}|x-y|^{2}
-\theta(|x-\overline{x}|^{4}+|y-\overline{x}|^{4})-\beta
(t-\overline{t})^2.
\end{equation}
By the linear growth assumption on $u$ and $w$, there exists a
$(t_{\epsilon},x_{\epsilon},y_{\epsilon})\in [t_0,T]\times B_R
\times B_R $, for $R$ large enough, such that:
$$\Phi_{\epsilon}(t_{\epsilon},x_{\epsilon},y_{\epsilon})=\max\limits_{(t,x,y)}\Phi_{\epsilon}(t,x,y).$$
On the other hand, from
$2\Phi_{\epsilon}(t_{\epsilon},x_{\epsilon},y_{\epsilon})\geq
\Phi_{\epsilon}(t_{\epsilon},x_{\epsilon},x_{\epsilon})+\Phi_{\epsilon}(t_{\epsilon},y_{\epsilon},y_{\epsilon})$,
we have
\begin{equation}
\frac{1}{\epsilon}|x_{\epsilon} -y_{\epsilon}|^{2} \leq
(1-\lambda)(u(t_{\epsilon},x_{\epsilon})-u(t_{\epsilon},y_{\epsilon}))+(w(t_{\epsilon},x_{\epsilon})-w(t_{\epsilon},y_{\epsilon})),
\end{equation}
and consequently $\frac{1}{\epsilon}|x_{\epsilon}
-y_{\epsilon}|^{2}$ is bounded, and as $\epsilon\rightarrow 0$,
$|x_{\epsilon} -y_{\epsilon}|\rightarrow 0$. Since $u$ and $w$ are
uniformly continuous on $[0,T]\times \overline{B}_R$, then
$\frac{1}{2\epsilon}|x_{\epsilon} -y_{\epsilon}|^{2}\rightarrow 0$
as
$\epsilon\rightarrow 0.$\\
Since
 $$(1-\lambda)u(\overline{t},\overline{x})-w(\overline{t},\overline{x}) \leq
 \Phi_{\epsilon}(t_{\epsilon}x_{\epsilon},y_{\epsilon})\leq (1-\lambda)u(t_{\epsilon},x_\epsilon)-w(t_{\epsilon},y_\epsilon),$$
it follow as $\lambda\rightarrow 0$ and the continuity of $u$ and
$w$ that, up to a subsequence,
 \begin{equation}\label{subsequence}
 (t_\epsilon,x_\epsilon,y_\epsilon)\rightarrow (\overline{t},\overline{x},\overline{x}).
 \end{equation}
Next let us show that $t_{\epsilon} <T.$ Actually if $t_{\epsilon}
=T$ then,
$$
\Phi_{\epsilon}(\overline{t},\overline{x},\overline{x})\leq
\Phi_{\epsilon}(T,x_{\epsilon},y_{\epsilon}),$$ and,
$$
(1-\lambda)u(\overline{t},\overline{x})-w(\overline{t},\overline{x})\leq
(1-\lambda)\exp(T)G_1(x_{\epsilon}) -\exp(T)G_1(y_{\epsilon})- \beta
(T-t_{\epsilon})^2,
$$
since $u(T,x_{\epsilon})=\exp(T)G_1(x_{\epsilon})$,
$w(T,y_{\epsilon})=\exp(T)G_1(y_{\epsilon})$ and $G_1$ is uniformly
continuous on $\overline{B}_R$. Then as $\lambda\rightarrow 0$ we
have,
$$
\begin{array}{ll}
\eta &\leq - \beta (T-\overline{t})^2\\ \eta &< 0,
\end{array}
$$
which yields a contradiction and we have $t_\epsilon \in [t_0,T)$.
We now claim that:
\begin{equation}
\label{visco-comp11}w(t_\epsilon,y_\epsilon)-\inf\limits_{\xi\in
E}[w(t_\epsilon,y_\epsilon+g(t_\epsilon,y_\epsilon,\xi))+\exp(t_\epsilon)C(t_\epsilon,y_\epsilon,\xi)]
< 0.
\end{equation}
Indeed if
$$w(t_\epsilon,y_\epsilon)-\inf\limits_{\xi\in
E}[w(t_\epsilon,y_\epsilon+g(t_\epsilon,y_\epsilon,\xi))+\exp(t_\epsilon)C(t_\epsilon,y_\epsilon,\xi)]\geq
0,$$  from the subsolution property of $u(t_\epsilon,x_\epsilon)$,
we have
$$ u(t_\epsilon,x_\epsilon)-\inf\limits_{\xi\in
E}[u(t_\epsilon,x_\epsilon+g(t_\epsilon,x_\epsilon,\xi))+\exp(t_\epsilon)C(t_\epsilon,x_\epsilon,\xi)]\leq
0,
$$
then there exists $\xi_1\in E$ such that:
$$w(t_\epsilon,y_\epsilon)-w(t_\epsilon,y_\epsilon+g(t_\epsilon,y_\epsilon,\xi_1))-\exp(t_\epsilon)C(t_\epsilon,y_\epsilon,\xi_1)\geq
0,$$
$$u(t_\epsilon,x_\epsilon)-u(t_\epsilon,x_\epsilon+g(t_\epsilon,x_\epsilon,\xi_1))-\exp(t_\epsilon)C(t_\epsilon,x_\epsilon,\xi_1)\leq0.$$It
follows that:
$$\begin{array}{ll}(1-\lambda)(u(t_\epsilon,x_\epsilon)- w(t_\epsilon,y_\epsilon) -[(1-\lambda)u(t_\epsilon,x_\epsilon+g(t_\epsilon,x_\epsilon,\xi_1))
-w(t_\epsilon,y_\epsilon+g(t_\epsilon,y_\epsilon,\xi_1))] \\\leq
(1-\lambda)\exp(t_\epsilon)C(t_\epsilon,x_\epsilon,\xi_1)-
\exp(t_\epsilon)C(t_\epsilon,y_\epsilon,\xi_1) .\end{array}$$ Now
since $C\geq \alpha >0$, then
$$\begin{array}{lll}(1-\lambda)(u(t_\epsilon,x_\epsilon)- w(t_\epsilon,y_\epsilon) -[(1-\lambda)u(t_\epsilon,x_\epsilon+g(t_\epsilon,x_\epsilon,\xi_1))
-w(t_\epsilon,y_\epsilon+g(t_\epsilon,y_\epsilon,\xi_1))]
\\<-\lambda\alpha+\exp(t_\epsilon)C(t_\epsilon,x_\epsilon,\xi_1)-\exp(t_\epsilon)C(t_\epsilon,y_\epsilon,\xi_1)
.\end{array}$$  But this contradicts the definition of
(\ref{comp-equ1}), since $C$, $u$ and $w$ is uniformly continuous on
$[0,T]\times \overline{B}_R$ and the claim (\ref{visco-comp11})
holds.

Next let us denote
\begin{equation}
\varphi_{\epsilon}(t,x,y)=\frac{1}{2\epsilon}|x-y|^{2}
+\theta(|x-\overline{x}|^{4}+|y-\overline{x}|^{4})+\beta
(t-\overline{t})^2.
\end{equation}
Then we have: \be \left\{
\begin{array}{lllll}\label{derive}
D_{t}\varphi_{\epsilon}(t,x,y)=2\beta(t-\overline{t}),\\
c+d=2\beta(t-\overline{t}),\\ D_{x}\varphi_{\epsilon}(t,x,y)=
\frac{1}{\epsilon}(x-y) +4\theta
(x-\overline{x})|x-\overline{x}|^{2}, \\
D_{y}\varphi_{\epsilon}(t,x,y)= -\frac{1}{\epsilon}(x-y) +
4\theta(y-\overline{x})|y-\overline{x}|^{2}.
\end{array}
\right. \ee Taking now into account (\ref{visco-comp11}), and the
definition of viscosity solution, we get:
\begin{equation}\begin{array}{lll}\label{vis_sub1}\displaystyle\min\limits_{\tau}[-c+(1-\lambda)u(t_\epsilon,x_\epsilon)
-\langle\frac{1}{\epsilon}(x_\epsilon-y_\epsilon) +4\theta
(x_\epsilon-\overline{x})|x_\epsilon-\overline{x}|^{2},\\\qquad\qquad\qquad\qquad
f(t_\epsilon,x_\epsilon,\tau)\rangle-(1-\lambda)\exp(t_\epsilon)\psi(t_\epsilon,x_\epsilon,\tau)]\leq0
\end{array}\end{equation}

 and

\begin{equation}\begin{array}{l}\label{vis_sub11}\displaystyle\min\limits_{\tau}\left[d+w(t_\epsilon,y_\epsilon)-\langle
\frac{1}{\epsilon}(x_\epsilon-y_\epsilon) -4\theta
(y_\epsilon-\overline{x})|y_\epsilon-\overline{x}|^{2},
f(t_\epsilon,y_\epsilon,\tau)\rangle-\exp(t_\epsilon)\psi(t_\epsilon,y_\epsilon,\tau)\right]\geq0
\end{array}\end{equation}

which implies that:
\begin{equation}
\begin{array}{llllll}
\label{viscder}
&-c-d+(1-\lambda)u(t_\epsilon,x_\epsilon)-w(t_\epsilon,y_\epsilon)\\&
\qquad\qquad\leq \displaystyle\min\limits_{\tau}\left[-\langle
\frac{1}{\epsilon}(x_\epsilon-y_\epsilon) -4\theta
(y_\epsilon-\overline{x})|y_\epsilon-\overline{x}|^{2},
f(t_\epsilon,y_\epsilon,\tau)\rangle -\exp(t_\epsilon)\psi(t_\epsilon,y_\epsilon,\tau)\right]\\
&\qquad\qquad
-\displaystyle\min\limits_{\tau}\left[-\langle\frac{1}{\epsilon}(x_\epsilon-y_\epsilon)
+4\theta (x_\epsilon-\overline{x})|x_\epsilon-\overline{x}|^{2},
f(t_\epsilon,x_\epsilon,\tau)\rangle-(1-\lambda)\exp(t_\epsilon)\psi(t_\epsilon,x_\epsilon,\tau)\right],\\&
\qquad\qquad\leq
\displaystyle\sup\limits_{\tau}\left[\langle\frac{1}{\epsilon}(x_\epsilon-y_\epsilon)
+4\theta (x_\epsilon-\overline{x})|x_\epsilon-\overline{x}|^{2},
f(t_\epsilon,x_\epsilon,\tau)\rangle+(1-\lambda)\exp(t_\epsilon)\psi(t_\epsilon,x_\epsilon,\tau)\right]\\
&\qquad\qquad -\displaystyle\sup\limits_{\tau}\left[\langle
\frac{1}{\epsilon}(x_\epsilon-y_\epsilon) -4\theta
(y_\epsilon-\overline{x})|y_\epsilon-\overline{x}|^{2},
f(t_\epsilon,y_\epsilon,\tau)\rangle
+\exp(t_\epsilon)\psi(t_\epsilon,y_\epsilon,\tau)\right],\\&
\qquad\qquad\leq
\displaystyle\sup\limits_{\tau}[\langle\frac{1}{\epsilon}(x_\epsilon-y_\epsilon)
,
f(t_\epsilon,x_\epsilon,\tau)-f(t_\epsilon,y_\epsilon,\tau)\rangle\\&\qquad\qquad+
\langle 4\theta
(x_\epsilon-\overline{x})|x_\epsilon-\overline{x}|^{2},
f(t_\epsilon,x_\epsilon,\tau)\rangle + \langle4\theta
(y_\epsilon-\overline{x})|y_\epsilon-\overline{x}|^{2},
f(t_\epsilon,y_\epsilon,\tau)\rangle\\
&\qquad\qquad

+(1-\lambda)\exp(t_\epsilon)\psi(t_\epsilon,x_\epsilon,\tau)-\exp(t_\epsilon)\psi(t_\epsilon,y_\epsilon,\tau)].
\end{array}
\end{equation}

Now, from (\ref{regbs1}), we get:
$$
\langle\frac{1}{\epsilon}(x_\epsilon-y_\epsilon),f(t_\epsilon,x_\epsilon,\tau)-f(t_\epsilon,y_\epsilon,\tau)\rangle
\leq \frac{C}{\epsilon}|x_\epsilon - y_\epsilon|^{2}.$$ Next $$
\langle 4\theta
(x_\epsilon-\overline{x})|x_\epsilon-\overline{x}|^{2},
f(t_\epsilon,x_\epsilon,\tau)\rangle\leq 4C\theta
|x_\epsilon||x_\epsilon-\overline{x}|^{3},$$and finally,
$$
\langle 4\theta
(y_\epsilon-\overline{x})|y_\epsilon-\overline{x}|^{2},
f(t_\epsilon,y_\epsilon,\tau)\rangle\leq 4C\theta
|y_\epsilon||y_\epsilon-\overline{x}|^{3}.$$Taking in to account
$$c+d=2\beta(t_\epsilon-\overline{t}).$$ So that by plugging into (\ref{viscder}) and note that $\lambda
>0$ we obtain:
\begin{equation}
\begin{array}{llll}
\label{viscder}
&-2\beta(t_\epsilon-\overline{t})+(1-\lambda)u(t_\epsilon,x_\epsilon)-w(t_\epsilon,y_\epsilon)\\&
\qquad\qquad\qquad\qquad\leq
\displaystyle\frac{C}{\epsilon}|x_\epsilon -
y_\epsilon|^{2}\\&\qquad\qquad\qquad\qquad+ 4C\theta
|x_\epsilon||x_\epsilon-\overline{x}|^{3} + 4C\theta
|y_\epsilon||y_\epsilon-\overline{x}|^{3}\\
&\qquad\qquad\qquad\qquad

+\sup\limits_{\tau}(1-\lambda)\exp(t_\epsilon)\psi(t_\epsilon,x_\epsilon,\tau)-\exp(t_\epsilon)\psi(t_\epsilon,y_\epsilon,\tau)]
\end{array}
\end{equation}By sending $\epsilon\rightarrow0$, $\lambda\rightarrow0$, $\theta
\rightarrow0$ and taking into account of the continuity of $\psi$,
we obtain $\eta \leq 0$ which is a contradiction. The proof of
Theorem \ref{uni1} is now complete. $\Box$
\medskip

%
%
%

\end{document}